\documentclass[a4paper,12pt]{article}

\input{mymacros.sty}

\newcommand{\Deltaone}{\Delta_1}

\newcommand{\DeltadPone}{\Delta_4}

\newcommand{\Mprim}{M_{\mathrm{prim}}}

\newcommand{\Tv}{T^\vee}
\newcommand{\bTv}{\bT^\vee}

\newcommand{\Wone}{W_1}
\newcommand{\None}{N_1}
\newcommand{\Mone}{M_1}
\newcommand{\bTvone}{\bT^\vee_1}
\newcommand{\bTone}{\bT_1}
\newcommand{\Tvone}{T^\vee_1}

\newcommand{\rhoNat}{\rho_{\mathrm Nat}}

\newcommand{\be}{\bse}

\title{A note on dimer models and McKay quivers}
\author{Kazushi Ueda and Masahito Yamazaki}
\date{}
\begin{document}

\maketitle




%
%

\begin{abstract}
We give one formulation of an algorithm
of Hanany and Vegh \cite{Hanany-Vegh}
which takes a lattice polygon as an input
and produces a set of isoradial dimer models.
We study the case of lattice triangles in detail
and discuss the relation with coamoebas
following Feng, He, Kennaway and Vafa \cite{Feng-He-Kennaway-Vafa}.
\end{abstract}


%
%

\section{Introduction}
Dimer models are introduced
in 1930s as a statistical mechanical model
for adsorption of di-atomic molecules
on the surface of a crystal
\cite{Fowler-Rushbrooke}.
A dimer model in this sense is a graph
consisting of the set $N$ of nodes
and the set $E$ of edges,
together with a function
$\scE : E \to \bR$
which represents the energy of adsorption.
A {\em perfect matching} is a subset $D$ of $E$
such that for any node $n \in N$,
there exists a unique edge $e \in D$ adjacent to $n$,
and one is interested in the asymptotic behavior
of the partition function
$$
 Z = \sum_{D \text{ : a perfect matching}}
  \exp \lb - \beta \sum_{e \in D} \scE(e) \rb
$$
as the graph becomes large.
Many interesting statistical mechanical models
such as domino tilings,
Ising models
and plane partitions
can be realized as a dimer model
on a particular planar graph.
See e.g. \cite{Kenyon_IDM} for a review
on dimer models.

In 1961, Kasteleyn \cite{Kasteleyn}
introduced a method to represent the partition function
as the determinant of a weighted adjacency matrix,
called the Kasteleyn matrix.
This is used to great effect
by Kenyon, Okounkov and Sheffield
\cite{Kenyon-Okounkov-Sheffield}
to study dimer models on periodic bipartite graph.
The determinant of the Kasteleyn matrix
with respect to a particular weighting is called
the {\em characteristic polynomial},
and plays an essential role in their work.
It is a Laurent polynomial in two variables
whose Newton polygon,
called the {\em characteristic polygon},
contains the information
of possible asymptotic behavior of the periodic dimer model
on the universal cover of the torus.


Recent advances in string theory
have uncovered a new connection
between dimer models and geometry,
%
culminating in the proposal that
a bicolored isoradial graph on a torus produces
an AdS/CFT dual pair of a toric Sasaki-Einstein 5-manifold and
an $N = 1$ superconformal field theory in four dimensions
\cite{
Franco-Hanany-Martelli-Sparks-Vegh-Wecht_GTTGBT,
Franco-Hanany-Vegh-Wecht-Kennaway_BDQGT,
Hanany-Kennaway_DMTD}.
Isoradiality is needed for the existence of
{\em R-charges} in superconformal field theory,
which is dual to volumes of
Sasaki-Einstein manifolds
\cite{
Benvenuti-Pando_Zayas-Tachikawa,
Butti-Zaffaroni_RTD,
Butti-Zaffaroni_TGQGT,
Kato_ZFSFT,
Lee-Rey,
Martelli-Sparks-Yau_GDAM,
Martelli-Sparks-Yau_SEMVM}.

The toric Calabi-Yau 3-fold,
obtained as the metric cone over the Sasaki-Einstein 5-manifold,
is determined by the characteristic polygon.
The dual superconformal field theory is determined
by a quiver with potential,
which is obtained as the dual graph of the bicolored graph.
Isoradiality implies
\cite{Mozgovoy-Reineke,
Broomhead,
Davison,
Ishii-Ueda_09}
that the path algebra of this quiver with potential
is a Calabi-Yau algebra
in the sense of Ginzburg
\cite{Ginzburg_CYA},
and the toric Calabi-Yau 3-fold can be obtained
as the moduli space of representations of the quiver
\cite{
Franco-Vegh_MSGTDM, Ishii-Ueda_08}.

An inverse construction of a bicolored isoradial graph on a torus
from a convex lattice polygon
is studied by Hanany and Vegh
\cite{Hanany-Vegh}.
Following their work,
we formulate the {\em linear Hanany-Vegh algorithm},
which takes a convex lattice polygon as an input
and produces a set of bicolored graphs on a torus.
The adjective {\em linear} comes from the fact that
{\em zig-zag paths} in the resulting graph
behave like lines,
which enables one to show
the following:

\begin{theorem} \label{th:isoradial}
Let $\Delta$ be a convex lattice polygon and
$G$ be a bicolored graph obtained from $\Delta$
by the linear Hanany-Vegh algorithm.
Then $G$ is isoradial and
$\Delta$ coincides with
the characteristic polygon of $G$
up to translation.
\end{theorem}

Although it is difficult to enumerate the output
of the linear Hanany-Vegh algorithm in general,
one can give a complete description
when the input is a triangle:

\begin{theorem} \label{th:triangle-McKay}
If the input of the linear Hanany-Vegh algorithm
is a lattice triangle,
then the output
consists of a unique graph
which is dual to the McKay quiver.
\end{theorem}

See \cite{Stienstra_dessins} for another attempt
to give a mathematical formulation
of the algorithm of Hanany and Vegh,
and \cite{Gulotta, Ishii-Ueda_09, Stienstra_CPA}
for an alternative algorithm due to Gulotta.

An important aspect of dimer models is
their relation with mirror symmetry.
The mirror of the toric Calabi-Yau 3-fold is
an algebraic curve in $(\bCx)^2$,
which governs the {\em limit shape}
of the melting crystal model
associated with the bicolored graph
\cite{Okounkov-Reshetikhin-Vafa_QCYCC, Ooguri-Yamazaki_ECG}.
The {\em coamoeba} of this curve is defined
by Passare and Tsikh
as its image by the argument map
\begin{equation} \label{eq:arg_map}
\begin{array}{cccc}
 \Arg : & (\bCx)^2 & \to & (\bR / \bZ )^2 \\
 & \vin & & \vin \\
 & (x, y) & \mapsto &
  \dfrac{1}{2 \pi} \lb \arg(x), \arg(y) \rb,
\end{array}
\end{equation}
whose behavior is expected by Feng, He, Kennaway and Vafa
\cite{Feng-He-Kennaway-Vafa}
to be described by the bicolored graph.

We show the following in this paper:

\begin{theorem} \label{th:main}
Let $\Delta$ be a lattice triangle and
$W$ be the Laurent polynomial obtained
as the sum of monomials corresponding
to the vertices of $\Delta$.
Then the bicolored graph obtained
in Theorem \ref{th:triangle-McKay} is
a deformation retract of the coamoeba.
\end{theorem}

A more detailed description of the coamoeba in this case
in terms of the McKay quiver of an abelian subgroup of $\SL_3(\bC)$
is given in Theorem \ref{th:coamoeba_triangle}.
The proof is based on the identification of the torus $(\bCx)^2$ 
with the dual of a maximal torus of $\SL_3(\bC)$.
This point of view will be used
in \cite{Ueda-Yamazaki_toricdP}
to study equivariant homological mirror symmetry
for two-dimensional toric Fano stacks.

The organization of this paper is as follows:
In Section \ref{sc:def},
we recall basic definitions on dimer models.
In Section \ref{sc:algorithm},
we formulate the linear Hanany-Vegh algorithm.
Theorem \ref{th:isoradial} is proved
in Section \ref{sc:inverse}.
In Section \ref{sc:McKay},
we recall the identification of McKay quivers
with hexagonal tilings of a real 2-torus
and prove Theorem \ref{th:triangle-McKay}.
We discuss the asymptotic behavior of coamoebas
in Section \ref{sc:asymptotic},
and study the case of the sum of three monomials
in detail in Section \ref{sc:coamoeba}.

\section{Characteristic polygon of a dimer model} \label{sc:def}

Let $N = \bZ^2$ be a free abelian group of rank two and
$M = \Hom(N, \bZ)$ be the dual group.
Put
$
 M_\bR = M \otimes \bR
$
and
$
 \Tv = M_\bR / M.
$
We fix the standard orientation and
the Euclidean inner product on $M_\bR \cong \bR^2$.
\begin{itemize}
 \item
A {\em graph} on $\Tv$ consists of
\begin{itemize}
 \item
a finite subset $N \subset \Tv$
called the set of {\em nodes}, and
 \item
another finite set $E$
called the set of {\em edges},\\
consisting of continuous maps
$e : [0, 1] \to \Tv$,
\end{itemize}
such that
$$
 e(0) \in N, \quad
 e(1) \in N, \quad
 e((0,1)) \cap N = \emptyset,
$$
and
$$
 e((0, 1)) \cap e'((0, 1)) = \emptyset
$$
for not necessarily distinct
edges $e, e' \in E$.
 \item
A pair $(n, e)$ of a node and an edge of a graph
is said to be {\em adjacent} if
$e(0) = n$ or $e(1) = n$.
 \item
A graph is said to be {\em bipartite}
if the set of nodes
can be divided into two disjoint set
$N = B \sqcup W$
so that every edge is adjacent
both to an element of $B$ and $W$.
 \item
A {\em bicolored graph} is a bipartite graph
together with a choice of a division
$N = B \sqcup W$
satisfying the condition above.
The elements of $B$ and $W$ are said to be
{\em black} and {\em white} respectively.
 \item
A {\em face} of a graph on $\Tv$
is a connected component of the complement
of the set of (the images of) edges.
 \item
A bicolored graph on $\Tv$ is a {\em dimer model}
if every face is simply-connected.
 \item
A {\em perfect matching}
on a dimer model $G = (B, W, E)$
is a subset $D \subset E$
such that for every node $n \in B \cup W$,
there is a unique edge $e \in D$
adjacent to it.
\end{itemize}

Let $G$ be a dimer model and
consider the bicolored graph $\Gtilde$ on $M_\bR$
obtained from $G$ by pulling-back
by the natural projection $M_\bR \to \Tv$.
The set of perfect matchings of $G$ can naturally be identified
with the set of periodic perfect matchings of $\Gtilde$.
Fix a reference perfect matching $D_0$.
Then for any
perfect matching $D$,
the union $D \cup D_0$
divides $M_\bR$ into connected components.
The height function $h_{D, D_0}$ is
a locally-constant function on
$M_\bR \setminus (D \cup D_0)$
which increases (resp. decreases)
by $1$
when one crosses an edge $e \in D$
with the black (resp. white) vertex
on the right
or an edge $e \in D_0$
with the white (resp. black) vertex
on the right.
This rule determines the height function
up to an additive constant.
The height function may not be periodic
even if $D$ and $D_0$ are periodic,
and the {\em height change}
$h(D, D_0) = (h_x(D, D_0), h_y(D, D_0)) \in N$
of $D$ with respect to $D_0$
is defined as the difference
\begin{align*}
 h_x(D, D_0) &= h_{D, D_0}(p+(1,0)) - h_{D, D_0}(p), \\
 h_y(D, D_0) &= h_{D, D_0}(p+(0,1)) - h_{D, D_0}(p) 
\end{align*}
of the height function,
which does not depend on the choice of
$p \in M_\bR \setminus (D \cup D_0)$.
The dependence of the height change
on the choice of the reference matching
is given by
$$
 h(D, D_1) = h(D, D_0) - h(D_1, D_0)
$$
for any three perfect matchings $D$, $D_0$ and $D_1$.
We will often suppress the dependence of the height difference
on the reference matching
and just write $h(D) = h(D, D_0)$.
\begin{itemize}
 \item
The {\em characteristic polynomial} of $G$ is
the Laurent polynomial in two variables
defined by
$$
 Z(x, y)
  = \sum_{\text{$D$ : a perfect matching}}
     x^{h_x(D)} y^{h_y(D)}.
$$
 \item
The {\em characteristic polygon} is
the Newton polygon of the characteristic polynomial,
i.e. the convex hull of
$$
 \{ (h_x(D), h_y(D)) \in N
      \mid \text{$D$ is a perfect matching} \}.
$$
\end{itemize}

\begin{figure}[htbp]
\begin{tabular}{cc}

\begin{minipage}{.45 \linewidth}
\centering
\input{P2_graph.pst}
\caption{A dimer model}
\label{fg:P2_graph}
\end{minipage}

&

\begin{minipage}{.45 \linewidth}
\centering
\input{P2_matching_one.pst}
\caption{$D_1:(1,0)$}
\label{fg:P2_matching_one}
\end{minipage}

\\
& \\

\begin{minipage}{.45 \linewidth}
\centering
\input{P2_matching_two.pst}
\caption{$D_2:(0,1)$}
\label{fg:P2_matching_two}
\end{minipage}

&

\begin{minipage}{.45 \linewidth}
\centering
\input{P2_matching_three.pst}
\caption{$D_3:(-1,-1)$}
\label{fg:P2_matching_three}
\end{minipage}

\\
& \\

\begin{minipage}{.45 \linewidth}
\centering
\input{P2_matching_four.pst}
\caption{$D_4:(0,0)$}
\label{fg:P2_matching_four}
\end{minipage}

&

\begin{minipage}{.45 \linewidth}
\centering
\input{P2_matching_five.pst}
\caption{$D_5:(0,0)$}
\label{fg:P2_matching_five}
\end{minipage}

\\

& \\

\begin{minipage}{.45 \linewidth}
\centering
\input{P2_matching_six.pst}
\caption{$D_6:(0,0)$}
\label{fg:P2_matching_six}
\end{minipage}

&

\begin{minipage}{.45 \linewidth}
\centering
\input{P2_height_one.pst}
\caption{The height function for $D_1$
with respect to $D_4$}
\label{fg:P2_height_one}
\end{minipage}

\\
\end{tabular}
\end{figure}

\begin{figure}[htbp]
\begin{minipage}{\linewidth}
\centering
\input{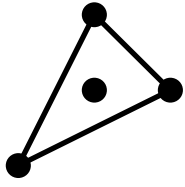}
\caption{The characteristic polygon}
\label{fg:Z3_diagram}
\end{minipage}
\end{figure}

As an example,
consider the dimer model
given in Figure \ref{fg:P2_graph}.
This dimer model has six perfect matchings
$D_1, \dots, D_6$
shown in Figures \ref{fg:P2_matching_one}--\ref{fg:P2_matching_six}.
The captions in these figures show the height changes
with respect to $D_4$,
and Figure \ref{fg:P2_height_one} shows
an example of a height function.
The characteristic polygon is shown in Figure \ref{fg:Z3_diagram}.

\section{Description of the algorithm}
 \label{sc:algorithm}
An element $m \in M$ is {\em primitive}
if it cannot be written as $k m'$ where $m' \in M$ and
$k$ is an integer greater than one.
The set of primitive elements of $M$ will be denoted by
$\Mprim$.
\begin{itemize}
 \item
An {\em oriented line} on $\Tv$ is
a pair $L = (\pi(\ell), m)$
of the image $\pi(\ell)$ of a line $\ell$ in $M_\bR$
by the natural projection
$
 \pi : M_\bR \to \Tv = M_\bR / M
$
and a primitive element $m \in \Mprim$
in the tangent space of $\ell$.
 \item
The map
$
 (\pi(\ell), m) \mapsto m
$
from the set of oriented lines
to $\Mprim$ is denoted by $p$.
 \item
The image $p(L)$ is called the {\em slope} of $L$.
 \item
The set of oriented lines on $\Tv$ can be identified
with the product $\Mprim \times S^1$,
where $\Mprim$ parametrizes the slope and
$S^1$ parametrizes the position of $\pi(\ell)$
in the direction perpendicular to the slope.
 \item
An {\em oriented line arrangement} on $\Tv$
is a finite set of oriented lines on $\Tv$.
 \item
An oriented line arrangement on $\Tv$ is {\em simple}
if no three lines intersect at one point.
 \item
A {\em cell} of an oriented line arrangement is
a connected component of the complement of the union of lines.
 \item
An {\em edge} of a cell
is a connected component of the intersection
of the boundary of the cell and
a line in the arrangement.
 \item
A {\em vertex} of a line arrangement is an intersection point
of two lines.
\end{itemize}
The following concepts are at the heart of the algorithm of
Hanany and Vegh:
\begin{itemize}
 \item
An edge of a cell of an oriented line arrangement has two orientations;
one comes from the orientation of the line,
and the other is the boundary orientation of the cell.
A cell is said to be {\em white} if these two orientations agree
for any of its edge,
and {\em black} if they are opposite
for any of its edge.
 \item
A cell is {\em colored} if it is either black or white.
Note that an edge in an oriented line arrangement bound
two cells, and at most one of them can be colored.
 \item
An oriented line arrangement is {\em admissible}
if every edge bounds a colored cell.
\end{itemize}
Let $\Delta \subset N_\bR$ be a convex lattice polygon,
i.e. the convex hull of a finite subset of $N$.
\begin{itemize}
 \item
An {\em edge} of $\Delta$
is a connected component of
$\partial \Delta \setminus (\partial \Delta \cap N)$.
The set of edges of $\Delta$ will be denoted by $\scE$.
 \item
The primitive outward normal vector to an edge $e$
of $\Delta$ is denoted by $m(e) \in \Mprim$.
 \item
The map $\Delta \mapsto \lc m(e) \rc_{e \in \scE}$ is a bijection
from the set of convex lattice polygons up to translations
to the set of finite collections of elements of $\Mprim$
summing up to zero.
 \item
An oriented line arrangement $\scA$ is said to be
{\em associated} with $\Delta$
if $\Delta$ corresponds to $\lc p(L) \rc_{L \in \scA}$
by the bijection above.
\end{itemize}
For example,
the triangle $\Delta_2$
shown in Figure \ref{fg:delta_two}
has four edges.
There are two combinatorially distinct ways
shown in Figure \ref{fg:delta_two_lines_one}
and Figure \ref{fg:delta_two_lines_two}
to arrange four lines
in the directions of outward normal vectors
of the edges of $\Delta_2$.
In Figure \ref{fg:delta_two_lines_one},
all edges bound colored polygons,
whereas in Figure \ref{fg:delta_two_lines_two},
the edges drawn in dotted lines
do not bound any colored polygon.
Hence the polygon $\Delta_2$ in Figure \ref{fg:delta_two}
has a unique admissible arrangement.

\begin{figure}[htbp]

\begin{tabular}{cc}


\multicolumn{2}{c}{
\begin{minipage}{.90\linewidth}
\centering
\input{delta_two.pst}
\caption{A triangle $\Delta_2$}
\label{fg:delta_two}
\end{minipage}
}
\\

& \\

\begin{minipage}{.45 \linewidth}
\centering
\input{delta_two_lines_one.pst}
\caption{The admissible arrangement}
\label{fg:delta_two_lines_one}
\end{minipage}

&

\begin{minipage}{.45\linewidth}
\centering
\input{delta_two_lines_two.pst}
\caption{The non-admissible arrangement}
\label{fg:delta_two_lines_two}
\end{minipage}

\\
& \\

\multicolumn{2}{c}{
\begin{minipage}{.90\linewidth}
\centering
\input{delta_two_graph.pst}
\caption{The dimer model}
\label{fg:delta_two_graph}
\end{minipage}
}

\\
\end{tabular}

\end{figure}

Now we define the map
$\frakD$ from the set of admissible arrangements
to the set of dimer models:
\begin{definition}
Let $\scA$ be an admissible oriented line arrangement.
Then the dimer model $\frakD(\scA) = (B, W, E)$ associated with $\scA$
is defined as follows:
\begin{itemize}
 \item
The set $W$ of white nodes
is the set of the centers of gravities of white cells.
 \item
The set $B$ of black nodes
is the set of the centers of gravities of black cells.
 \item
Two nodes are connected by a straight line segment
if the corresponding cells share a vertex.
\end{itemize}
\end{definition}

As an example,
Figure \ref{fg:delta_two_graph} shows the dimer model
associated with the admissible arrangement
in Figure \ref{fg:delta_two_lines_one}.

Now we can formulate the linear Hanany-Vegh algorithm:

\begin{definition}
The {\em linear Hanany-Vegh algorithm} is defined as follows:
\begin{itemize}
 \item
Take a convex lattice polygon $\Delta \subset N_\bR$ as an input.
 \item
The set $\frakA(\Delta)$ of oriented line arrangements $\{ L_a \}_a$
associated with $\Delta$
is a real $n$-torus,
where $n$ is the number of edges of $\Delta$.
 \item
The subset of $\frakA(\Delta)$ corresponding to simple arrangements
will be denoted by $\frakU(\Delta)$.
The complement $\frakA(\Delta) \setminus \frakU(\Delta)$
is a closed subset of real codimension one.
 \item
Take a representative
$
 \frakV(\Delta) \subset \frakU(\Delta)
$
of $\pi_0(\frakU(\Delta))$ and let
$
 \frakW(\Delta) \subset \frakV(\Delta)
$
be the subset consisting of admissible arrangements.
 \item
The output is the set $\frakD(\frakW(\Delta))$ of dimer models
associated with arrangements in $\frakW(\Delta)$.
\end{itemize}
\end{definition}

As an example,
take the convex hull $\Delta_4$ of
$$
 v_1 = (1,0), \ 
 v_2 = (0,1), \ 
 v_3 = (-1,0), \ \text{and} \ 
 v_4 = (-1, -1)
$$
shown in Figure \ref{fg:dP1_diagram}.
The primitive outward normal vectors
of the edges of $\DeltadPone$ are given by
$$
 m_1 = (1,1), \ 
 m_2 = (-1,1), \ 
 m_3 = (-1,0), \ \text{and} \ 
 m_4 = (1, -2).
$$

\begin{figure}[htbp]
\begin{tabular}{cc}
\begin{minipage}{.45 \linewidth}
\centering
\input{dP1_diagram.pst}
\caption{The lattice polygon $\DeltadPone$
with its primitive normal vectors}
\label{fg:dP1_diagram}
\end{minipage}

&

\begin{minipage}{.45 \linewidth}
\centering
\input{dP1_lines.pst}
\caption{The unique admissible arrangement
corresponding to $\DeltadPone$}
\label{fg:dP1_lines}
\end{minipage}
\end{tabular}
\end{figure}

\begin{lemma} \label{lem:unique_dP1}
There is a unique admissible oriented line arrangement
associated with $\DeltadPone$.
\end{lemma}
\begin{proof}
Figure \ref{fg:three_lines} shows
the unique arrangement of three oriented lines
on the torus in the directions of $m_1$, $m_2$, and $m_4$.
There are five combinatorially distinct ways
to insert a line into Figure \ref{fg:three_lines}
in the direction of $m_3$
shown in Figure \ref{fg:four_lines_one}.
\begin{figure}
\begin{tabular}{cc}
\begin{minipage}{.45 \linewidth}
\centering
\input{three_lines.pst}
\caption{The unique arrangement of three lines
in the directions of $m_1$, $m_2$ and $m_4$}
\label{fg:three_lines}
\end{minipage}
&
\begin{minipage}{.45 \linewidth}
\centering
\input{four_lines_one.pst}
\caption{Five ways to insert the fourth line
in the direction of $m_3$}
\label{fg:four_lines_one}
\end{minipage}
\end{tabular}
\end{figure}
In addition,
since two of six intersection points
of the three lines
in Figure \ref{fg:three_lines}
are on the same horizontal level,
there are two combinatorially distinct ways
shown in Figure \ref{fg:four_lines_two}
and Figure \ref{fg:four_lines_three}
to perturb them a little and
insert the fourth line.
\begin{figure}[htbp]
\begin{tabular}{cc}
\begin{minipage}{.45 \linewidth}
\centering
\input{four_lines_two.pst}
\caption{one way to perturb three lines
to insert the fourth line}
\label{fg:four_lines_two}
\end{minipage}
&
\begin{minipage}{.45 \linewidth}
\centering
\input{four_lines_three.pst}
\caption{another way to perturb three lines
to insert the fourth line}
\label{fg:four_lines_three}
\end{minipage}
\end{tabular}
\end{figure}
It is easy to see that
out of these seven arrangements,
only the one shown in Figure \ref{fg:dP1_lines}
is admissible.
\end{proof}

\begin{figure}[htbp]
\centering
\input{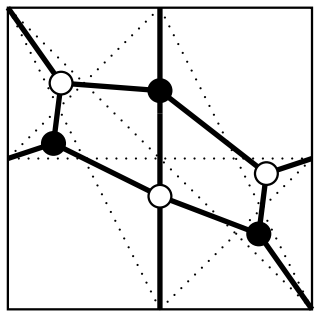}
\caption{The dimer model}
\label{fg:dP1_graph}
\end{figure}

Hence the output of the linear Hanany-Vegh algorithm consists of
only one dimer model shown in Figure \ref{fg:dP1_graph}
in this case.

\section{Zig-zag paths and isoradiality}
 \label{sc:inverse} 

The following notion is due to
Duffin \cite{Duffin} and Mercat \cite{Mercat_DRSIM}:

\begin{definition} \label{df:isoradial}
A dimer model is {\em isoradial}
if one can choose an embedding of the graph
into the torus $\Tv$
so that every face of the graph is a polygon
inscribed in a circle of a fixed radius
with respect to a flat metric on $\Tv$. 
Here, the circumcenter of any face must be contained
in the face.
\end{definition}

The following notion is introduced
by Kenyon and Schlenker \cite{Kenyon-Schlenker}.

\begin{definition} \label{df:zigzag}
A {\em zig-zag path} is a path on a dimer model
which makes a maximum turn to the right on a white node
and maximum turn to the left on a black node.
\end{definition}

\begin{figure}[htbp]
\centering
\input{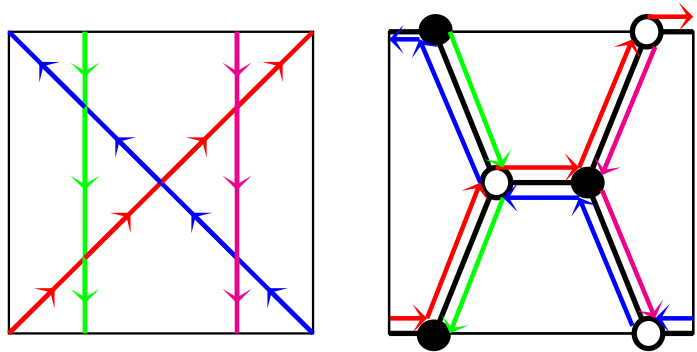}
\caption{An admissible arrangement and zig-zag paths}
\label{fg:delta_two_zigzag}
\end{figure}

If a dimer model comes from a lattice polygon
through the linear Hanany-Vegh algorithm,
then the set of zig-zag path can naturally be identified
with the set of oriented lines in the admissible arrangement.
Figure \ref{fg:delta_two_zigzag} shows
an example of this identification.

A dimer model is isoradial
if and only if zig-zag paths behave like straight lines:

\begin{theorem}[{Kenyon and Schlenker \cite[Theorem 5.1]{Kenyon-Schlenker}}]
 \label{th:Kenyon-Schlenker}
A dimer model is isoradial
if and only if the following conditions are satisfied:
\begin{enumerate}
 \item
Every zig-zag path is a simple closed curve.
 \item
The lift of any pair of zig-zag paths to the universal cover
of the torus intersect at most once.
\end{enumerate}
\end{theorem}

Since zig-zag paths in a dimer model
obtained through the linear Hanany-Vegh algorithm
clearly satisfies these conditions,
Theorem \ref{th:Kenyon-Schlenker} immediately implies the following:

\begin{corollary} \label{cor:isoradial}
A dimer model obtained by the linear Hanany-Vegh algorithm is isoradial.
\end{corollary}

The following notion is slightly weaker than
isoradiality:

\begin{definition}[{\cite[Definition 5.2]{Ishii-Ueda_09}}]
 \label{def:consistency}
A dimer model is {\em consistent} if
\begin{itemize}
 \item
no zig-zag path has a self-intersection on the universal cover,
 \item
no zig-zag path on the universal cover is a closed path, and
 \item
no pair of zig-zag paths intersect each other
on the universal cover in the same direction more than once.
\end{itemize}
\end{definition}


One can formulate the {\em consistent Hanany-Vegh algorithm}
by enlarging the domain of the brute-force search
from the set of oriented line arrangements
to the set of arrangements of oriented curves
satisfying the three conditions
in Definition \ref{def:consistency}.
The orientations are given by the homology classes of the curves
instead of the tangent spaces to the lines.
Note that for a convex lattice polygon $\Delta$,
the number of combinatorial types of arrangements of oriented curves
associated with $\Delta$
satisfying the conditions in Definition \ref{def:consistency}
is finite.
One can also formulate the {\em isoradial Hanany-Vegh algorithm}
by restricting to arrangements satisfying the conditions
in Theorem \ref{th:Kenyon-Schlenker}.

The following theorem shows that
Hanany-Vegh algorithms are inverse algorithms:

\begin{theorem}[{Gulotta \cite[Theorem 3.3]{Gulotta},
Ishii and Ueda \cite[Corollary 8.3]{Ishii-Ueda_09}}]
If the set of zig-zag paths
on a consistent dimer model $G$ is associated
with a lattice polygon $\Delta$,
then the characteristic polygon of $G$ coincides with $\Delta$
up to translation.
\end{theorem}

\section{Dimer models and McKay quivers}
 \label{sc:McKay}

We recall the description of the McKay quiver
for an abelian subgroup of $\SL_3(\bC)$
as the dual of a hexagonal tiling on a torus
\cite{Nakamura_HSAGO, Reid_MC} and
prove Theorem \ref{th:triangle-McKay}
in this section.

\subsection{Quivers with potentials}

A {\em quiver} consists of
\begin{itemize}
 \item a set $V$ of vertices,
 \item a set $A$ of arrows, and
 \item two maps $s, t: A \to V$ from $A$ to $V$.
\end{itemize}
For an arrow $a \in A$,
the vertices $s(a)$ and $t(a)$
are said to be the {\em source}
and the {\em target} of $a$
respectively.
A {\em path} on a quiver
is an ordered set of arrows
$(a_n, a_{n-1}, \dots, a_{1})$
such that $s(a_{i+1}) = t(a_i)$
for $i=1, \dots, n-1$.
We also allow for a path of length zero,
starting and ending at the same vertex.
The {\em path algebra} $\bC Q$
of a quiver $Q = (V, A, s, t)$
is the algebra
spanned by the set of paths
as a vector space,
and the multiplication is defined
by the concatenation of paths;
$$
 (b_m, \dots, b_1) \cdot (a_n, \dots, a_1)
  = \begin{cases}
     (b_m, \dots, b_1, a_n, \dots, a_1) & s(b_1) = t(a_n), \\
      0 & \text{otherwise}.
    \end{cases}
$$
A {\em quiver with relations}
is a pair of a quiver
and a two-sided ideal $\scI$
of its path algebra.


Let $p = (a_n, \dots, a_1)$ be
an oriented cycle on the quiver,
i.e. the class of a path satisfying $s(a_1) = t(a_n)$
up to a cyclic change of the starting point.
For an arrow $b$,
the {\em derivative} of $p$ by $b$ is defined by
$$
 \frac{\partial p}{\partial b}
  = \sum_{i=1}^n \delta_{a_i, b}
     (a_{i-1}, a_{i-2}, \dots, a_1, a_n, a_{n-1}, \dots, a_{i+1}),
$$
where
$$
 \delta_{a, b}
  = \begin{cases}
      1 & a = b, \\
      0 & \text{otherwise}.
    \end{cases}
$$
This derivation extends
to formal sums of oriented cycles by linearity.

A {\em potential} on a quiver is a formal sum 
$\Phi$ of oriented cycles,
which defines an ideal $\scI = (\partial \Phi)$
of relations
generated by the derivatives $\partial \Phi / \partial a$
for all the arrows $a$ of the quiver;
$$
 (\partial \Phi) 
   = \left(
      \frac{\partial \Phi}{\partial a}
     \right)_{a \in A}.
$$
The most basic example is the quiver
with $V=\{ v \}$ and $A = \{x, y, z\}$.
Since there is only one vertex,
the maps $s$ and $t$ must be constant.
Consider the potential
$$
 \Phi = y x z - x y z.
$$
The corresponding ideal is given by
\begin{align*}
 (\partial \Phi)
  &= \left( \frac{\partial \Phi}{\partial x},
       \frac{\partial \Phi}{\partial y},
       \frac{\partial \Phi}{\partial z}
     \right) \\
  &= ( y z - z y, z x - x z, x y - y x ).
\end{align*}
The resulting path algebra with relations
is just the polynomial ring in three variables;
$$
 \bC \langle x, y, z \rangle / (\partial \Phi)
  \cong \bC[x, y, z].
$$
See Ginzburg \cite{Ginzburg_CYA}
and references therein
for more on algebras
associated with quivers with potentials.

\subsection{McKay quivers}

Let $G$ be a finite subgroup of $\GL_n(\bC)$.
The {\em McKay quiver} of $G$ is a quiver
whose set of vertices is the set
of irreducible representations of $G$,
and the number $a_{\sigma \tau}$
of arrows from
$\sigma$ to $\tau$
is given by the multiplicity
of $\tau$ in the tensor product of $\sigma$
with the dual representation
$\rhoNat^\vee$ of the natural representation of $G$;
$$
 \sigma \otimes \rhoNat^\vee
 = \bigoplus_{\tau} \tau^{\otimes a_{\sigma \tau}}.
$$
A path of length two from $\sigma$ to $\tau$
corresponds to an element of
$$
 \Hom_G(\sigma \otimes \rhoNat^\vee \otimes \rhoNat^\vee, \tau).
$$
The McKay quiver comes with relations
generated by the kernel of the map
$$
 \sigma \otimes \rhoNat^\vee \otimes \rhoNat^\vee
  \to \sigma \otimes (\Sym^2 \rhoNat^\vee).
$$
When $G$ is an abelian subgroup $A$ of $\SL_3(\bC)$,
the McKay quiver has the following description:
Assume that $A$ is contained in the diagonal subgroup,
and let
$$
 \pi_i : A \to \bCx, \qquad i = 1, 2, 3,
$$
be the projections to the $i$-th diagonal component,
viewed as a one-dimensional representation of $A$.
The the set of vertices is
the set $\Irrep(A)$ of irreducible representations of $A$,
and for each vertex $\rho \in \Irrep(A)$,
there exist three arrows
$x_\rho$, $y_\rho$, and $z_\rho$,
ending at the vertices
$\rho \otimes \pi_1^{\vee}$,
$\rho \otimes \pi_2^{\vee}$, and
$\rho \otimes \pi_3^{\vee}$ respectively.
The relations come from the potential
\begin{equation} \label{eq:McKay_potential}
 \Phi = x z y - x y z,
\end{equation}
where
$x = \sum_\rho x_\rho$,
$y = \sum_\rho y_\rho$, and
$z = \sum_\rho z_\rho$.

\subsection{Toric Calabi-Yau 3-folds associated with lattice polygons}
Let $\Ntilde = N \oplus \bZ$ be a free abelian group
of rank three and
$\Mtilde = M \oplus \bZ$ be the dual group.
For a lattice polygon $\Delta \subset N_\bR$,
let $\sigma$ be the cone
over $\Delta \times \{ 1 \} \subset \Ntilde_\bR$ and
$\Sigma$ be the fan in $\Ntilde_\bR$
consisting of $\sigma$ and its faces.
Let $\{ \vtilde_i \}_{i=1}^r$ be the set of primitive generators
of one-dimensional cones of $\Sigma$, and
$$
 \phitilde : \bZ^r \to \Ntilde
$$
be the map which send the $i$-th standard coordinate vector $e_i \in \bZ^r$
to $\vtilde_i \in \Ntilde$.
Let
$$
 \bTtilde
  = \Ntilde \otimes \bCx
  = \Spec \bC[\Mtilde]
$$
be a three-dimensional torus.
The group
$$
 K = \Ker(\phitilde \otimes \bCx : (\bCx)^r \to \bTtilde) \subset (\bCx)^r
$$
naturally acts on $\bC^r$,
and the toric variety $X$ associated with $\Sigma$
is the quotient of $\bC^r$ by the action of $K$;
$$
 X = \bC^r / K.
$$
The monoid ring
$$
 R = \bC[\sigma^\vee]
$$
of the dual cone
$$
 \sigma^\vee
  = \lc m \in \Mtilde \mid
         \la m, n \ra \ge 0, \  n \in \sigma \rc
  \subset \Mtilde
$$
is the coordinate ring of $X$;
$$
 X = \Spec R.
$$
When
$\Delta$ is the convex hull of
\begin{align*}
 v_1 = (p, q), \qquad
 v_2 = (r, s), \qquad
 v_3 = (0, 0),
\end{align*}
then the map
$
 \phitilde : \bZ^3 \to \Ntilde
$
is given by
$$
 (a, b, c) \mapsto (p a + r b, q a + r b, a + b + c),
$$
so that the map
$
 \phitilde \otimes \bCx : \bZ^3 \to \bTtilde
$
is given by
$$
 (\alpha, \beta, \gamma)
  \mapsto (\alpha^p \beta^r, \alpha^q \beta^s, \alpha \beta \gamma).
$$
It follows that
the toric variety associated with $\Delta$
is the quotient
$$
 X = \bC^3 / A
$$
where
\begin{align*}
 A
 &= \lc (\alpha, \beta, \gamma) \in (\bCx)^3 \mid
      \alpha^p \beta^r = \alpha^q \beta^s = \alpha \beta \gamma = 1 \rc.
\end{align*}

\subsection{McKay quivers on the torus}
The McKay quiver of an abelian subgroup of $\SL_3(\bC)$
can naturally be drawn on a real 2-torus as follows:
Let $\None = \bZ^2$ be a free abelian group of rank two and
$\Mone = \Hom(\None, \bZ)$ be the dual group.
Identify the torus
$$
 \bTone = \None \otimes \bCx = \Spec \bC[\Mone]
$$
with a maximal torus
$$
 \{ \diag(\alpha, \beta, \gamma) \in \SL_3(\bC)
  \mid \alpha \beta \gamma = 1 \}
$$
of $\SL_3(\bC)$,
so that
the group
$$
 A = \lc \diag(\alpha, \beta, \gamma) \in \SL_3(\bC) \mid
      \alpha^p \beta^r = \alpha^q \beta^s = \alpha \beta \gamma = 1 \rc
$$
is identified with the kernel of the map
$$
\begin{array}{cccc}
 \phi \otimes \bCx : & \bTone & \to & \bT \\
 & \vin & & \vin \\
 & (\alpha, \beta) & \mapsto &
  (\alpha^p \beta^r, \alpha^q \beta^s),
\end{array}
$$
where $\phi : \None \to N$ is the linear map
given by the matrix
$$
 P =
 \begin{pmatrix}
  p & r \\
  q & s
 \end{pmatrix}.
$$
The short exact sequence
$$
 1 \to A \to \bTone \xto{\phi \otimes \bCx} \bT \to 1
$$
of abelian groups induces an exact sequence
$$
 0 \to M \xto{\psi} \Mone \to \Irrep(A) \to 0
$$
of characters,
where $\psi$ is the adjoint map of $\phi : \None \to N$
represented by the transposed matrix of $P$.

\begin{figure}[htbp]
\centering
\input{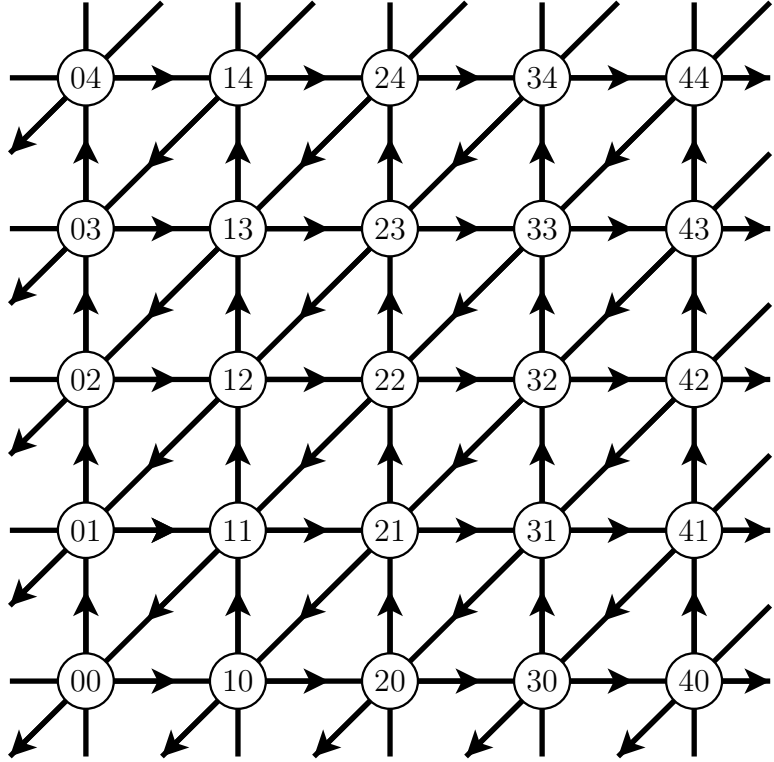}
\caption{The periodic quiver}
\label{fg:periodic_quiver}
\end{figure}

Now consider the infinite quiver
in Figure \ref{fg:periodic_quiver}
drawn on $\Mone \otimes \bR \cong M_\bR$,
whose set of vertices is $\Mone$ and
whose set of arrows consists of
$x_{i, j}$, $y_{i, j}$ and $z_{i, j}$
for $(i, j) \in \Mone$ such that
$$
 s(x_{i, j}) = s(y_{i, j}) = s(z_{i, j}) = (i, j)
$$
and
$$
 t(x_{i, j}) = (i+1, j), \quad
 t(y_{i, j}) = (i, j+1), \quad
 t(z_{i, j}) = (i-1, j-1).
$$
The McKay quiver for $A$ is the quotient of this quiver
by the natural action of $M$.

\begin{figure}[htbp]
\centering
\input{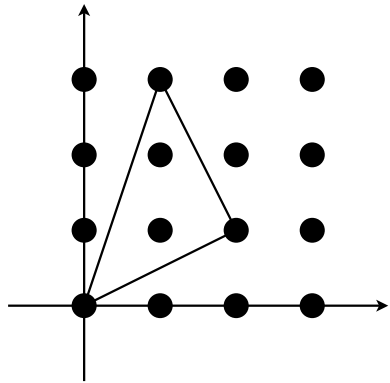}
\caption{The triangle $\Delta_5$}
\label{fg:Z5_diagram}
\end{figure}

\begin{figure}[htbp]
\centering
\input{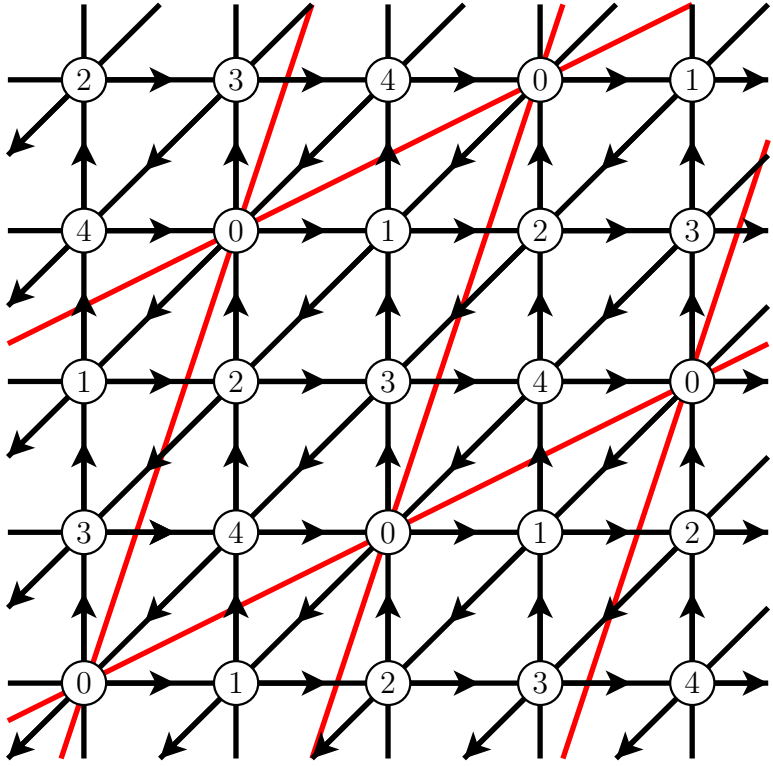}
\caption{The McKay quiver}
\label{fg:Z5_quiver}
\end{figure}

As an example,
consider the case when
$\Delta$ is the convex hull of $(2, 1)$, $(1, 3)$ and $(0, 0)$
shown in Figure \ref{fg:Z5_diagram}.
The map $\phi$ is presented by the matrix
$$
 P =
 \begin{pmatrix}
  2 & 1 \\
  1 & 3
 \end{pmatrix},
$$
and the subgroup $A$ is given by
$$
 A
  = \lc (\alpha, \beta, \gamma) \in \bTone
     \mid \alpha^2 \beta = \alpha \beta^3
           = \alpha \beta \gamma = 1 \rc,
$$
which is isomorphic to the image of
$$
\begin{array}{cccc}
 \rhoNat : & \bZ / 5 \bZ & \to & \bT \\
 & \vin & & \vin \\
 & [1] & \mapsto & (\zeta, \zeta^3, \zeta)
\end{array}
$$
where
$
 \zeta = \exp(2 \pi \sqrt{-1} / 5).
$
The quotient of the periodic quiver
in Figure \ref{fg:periodic_quiver}
by the action of $M$
is shown in Figure \ref{fg:Z5_quiver},
where fundamental regions of the action
are shown in red.

\subsection{A quiver with potential from a dimer model}

\begin{figure}[htbp]
\centering
\input{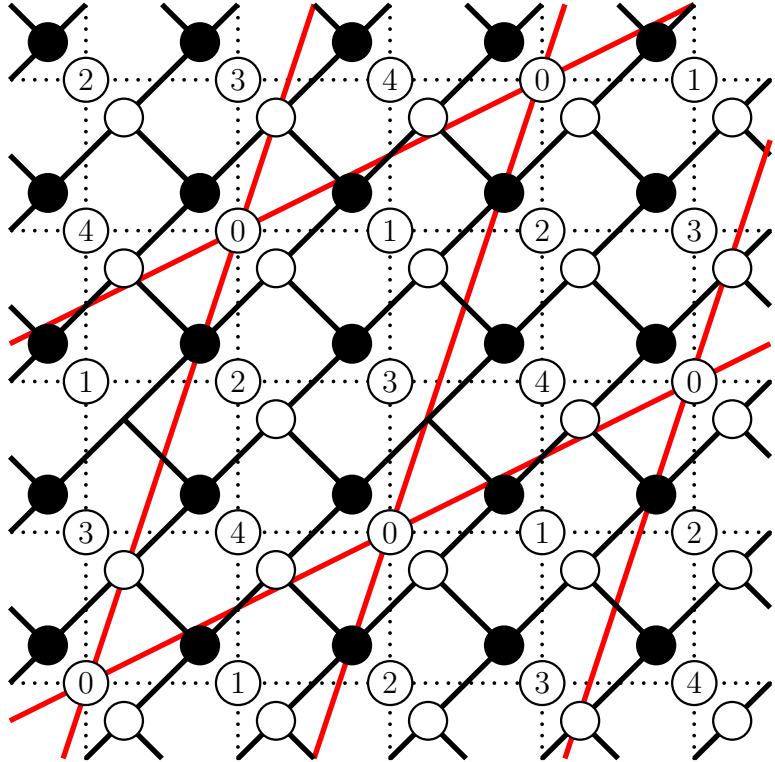}
\caption{The dual honeycomb graph}
\label{fg:Z5_graph}
\end{figure}

Note that the dual graph of the quiver
in Figure \ref{fg:Z5_quiver} is
a hexagonal tiling of the torus
$\Tv = M_\bR / M$
shown in Figure \ref{fg:Z5_graph}.
The colors of the nodes of the dual graph
is chosen so that
a white node is always on the right of an arrow.
Observe that the potential \eqref{eq:McKay_potential},
which gives the relations of the McKay quiver,
is the signed sum of cycles around the nodes of this graph;
\begin{equation} \label{eq:potential}
 \Phi = \sum_{w \in W} p_w - \sum_{b \in B} p_b.
\end{equation}
This motivates the following definition
of the quiver with potential associated with a dimer model
$G = (B, W, E)$:
\begin{itemize}
 \item
The set $V$ of vertices
is the set of faces of $G$.
 \item
The set $A$ of arrows
is the set $E$ of edges of the graph.
The directions of the arrows are determined
by the colors of the vertices of the graph,
so that the white vertex $w \in W$ is on the right
of the arrow.
 \item
A node $n \in B \coprod W$ determines
an oriented cycle $p_n$ of the quiver around it,
and the potential is defined
as the signed sum \eqref{eq:potential} over such cycles.
\end{itemize}
A hexagonal tiling of a torus
can naturally be identified with the McKay quiver
of an abelian subgroup of $\SL_3(\bC)$ in this way.

\subsection{The linear Hanany-Vegh algorithm on triangles}

Now we prove Theorem \ref{th:triangle-McKay}.
First note that
given three slopes,
there is only one way to form a polygon
whose edges has one of these slopes
and the boundary is positively oriented.
This unique arrangement is the triangle
shown in Figure \ref{fg:triangle_proof_one},
and an example of a polygon
with unoriented boundary is shown
in Figure \ref{fg:triangle_proof_six}.

\begin{figure}[htbp]
\begin{tabular}{cc}
\begin{minipage}{.45 \linewidth}
\centering
\input{triangle_proof_one.pst}
\caption{The unique positively-oriented polygon}
\label{fg:triangle_proof_one}
\end{minipage}
&
\begin{minipage}{.45 \linewidth}
\centering
\input{triangle_proof_six.pst}
\caption{an example of an unoriented pentagon}
\label{fg:triangle_proof_six}
\end{minipage}
\\
\ \\
\begin{minipage}{.45 \linewidth}
\centering
\input{triangle_proof_two.pst}
\caption{Two of the adjacent triangles}
\label{fg:triangle_proof_two}
\end{minipage}
&
\begin{minipage}{.45 \linewidth}
\centering
\input{triangle_proof_three.pst}
\caption{Two more triangles}
\label{fg:triangle_proof_three}
\end{minipage}
\\
\ \\
\begin{minipage}{.45 \linewidth}
\centering
\input{triangle_proof_four.pst}
\caption{Yet another triangle}
\label{fg:triangle_proof_four}
\end{minipage}
&
\begin{minipage}{.45 \linewidth}
\centering
\input{triangle_proof_five.pst}
\caption{The hexagonal tiling}
\label{fg:triangle_proof_five}
\end{minipage}
\end{tabular}
\end{figure}

For an arrangement of oriented lines to be admissible,
each of the vertices of the triangle
in Figure \ref{fg:triangle_proof_one}
must be contained in another triangle,
whose boundary is negatively oriented.
Again, there is a unique such arrangement,
and Figure \ref{fg:triangle_proof_two}
shows two of these adjacent triangles.

By repeating the same argument,
one obtains two new triangles
in Figure \ref{fg:triangle_proof_three},
and then another triangle
in Figure \ref{fg:triangle_proof_four}.
It follows that the hexagonal tiling of the the torus
which locally looks as
in Figure \ref{fg:triangle_proof_five}
is the unique output of the linear Hanany-Vegh algorithm
on a triangle $\Delta$.
It is clear that this tiling is dual
to the McKay quiver associated with $\Delta$,
and Theorem \ref{th:triangle-McKay} is proved.

\section{Coamoebas and Newton polygons}
\label{sc:asymptotic}

For a Laurent polynomial
$$
 W(x, y) = \sum_{(i,j) \in N} a_{ij} x^i y^j
$$
in two variables,
its {\em Newton polygon}
is defined as the convex hull
of $(i,j) \in N$ such that $a_{ij} \ne 0$;
$$
 \Delta = \Conv
 \{ (i,j) \in N
     \mid a_{ij} \neq 0 \}
 \subset N_\bR.
$$
It gives a regular map
$$
 W : \bTv \to \bC
$$
from the torus
$
 \bTv = M \otimes \bCx = \Spec \bC[N]
$
dual to $\bT = N \otimes \bCx$.
The {\em coamoeba} of $W^{-1}(0) \subset \bTv$
is its image under the argument map
$$
\begin{array}{cccc}
 \Arg : & \bTv & \to & \Tv \\
  & \vin & & \vin \\
  & (x, y) & \mapsto & \dfrac{1}{2 \pi}(\arg x, \arg y).
\end{array}
$$
For an edge $e$ of $\Delta$,
let $(n(e), m(e)) \in M$
be the primitive outward normal vector of $e$ and
$l(e)$ be the integer such that
the defining equation for the edge $e$ is given by
$$
 n(e) i + m(e) j = l(e).
$$
The {\em leading term} of $W$
with respect to the edge $e$ is defined by
$$
 W_e (x, y)
  = \sum_{n(e) i + m(e) j = l(e)} a_{ij} x^i y^j.
$$
This is indeed the leading term
if we put
$$
 (x, y)
   = (r^{n(e)} u, r^{m(e)} v),
 \qquad r \in \bR \ \text{and} \ 
 u, v \in \bCx
$$
and take the $r \to \infty$ limit;
$$
 P(r^{n(e)} u, r^{m(e)} v)
  = r^{l(e)} W_e(u, v) + O(r^{l(e)-1}).
$$
Now assume that for an edge $e$,
the leading term $W_e(x,y)$ is a binomial
$$
 W_e(x, y)
  = a_1 x^{i_1} y^{j_1}
     + a_2 x^{i_2} y^{j_2},
$$
where $a_1, a_2 \in \bC$ and
$(i_1, j_1), (i_2, j_2) \in N$.
Put $\alpha_i = \arg(a_i)$ for $i = 1, 2$ and
$$
 (x, y) =
   (r^{n(e)} |a_2| \be(\theta), r^{m(e)} |a_1| \be(\phi)),
$$
where $\be(x) = \exp(2 \pi \sqrt{-1} x)$
for $x \in \bR / \bZ$.
Then the leading behavior of $W$ as $r \to \infty$
is given by
\begin{equation} \label{eq:coamoeba_boundary}
 r^{l(e)} W_e(\be(\theta), \be(\phi))
  =  r^{l(e)} |a_1 a_2| \{ \be(\alpha_1 + i_1 \theta + j_1 \phi)
                + \be(\alpha_2 + i_2 \theta + j_2 \phi) \}.
\end{equation}
Hence the coamoeba of $W^{-1}(0)$ asymptotes in this limit
to the line
$$
 (\alpha_2 - \alpha_1) + (i_2 - i_1) \theta + (j_2 - j_1) \phi
  + \frac{1}{2} = 0 \mod \bZ
$$
on the torus $\Tv$.
This line will be called
an {\em asymptotic boundary} of the coamoeba of $W^{-1}(0)$.

The asymptotic boundary has a natural orientation
coming from the outward normal vector of the edge of $\Delta$.
To understand the role of this orientation,
take a pair of adjacent edges $e$ and $e'$ of $\Delta$
as in Figure \ref{fg:asympbone}
and consider the behavior of the coamoeba of $W^{-1}(0)$
near the intersection of asymptotic boundaries
corresponding to $e$ and $e'$.
Assume that the leading terms
corresponding to $e$ and $e'$ are binomials
\begin{eqnarray*}
  W_e(x, y)
   & = \be(\alpha) x^{i_1} y^{j_1}
        + \be(\beta) x^{i_2} y^{j_2}, \\
  W_{e'}(x, y)
   & = \be(\beta) x^{i_2} y^{j_2}
        + \be(\gamma) x^{i_3} y^{j_3}
\end{eqnarray*}
for some $\alpha, \beta, \gamma \in \bR/\bZ$
and $(i_1, j_1), (i_2, j_2), (i_3, j_3) \in N$.
Put
$$
 W_{e e'}(x, y) = \be(\alpha) x^{i_1} y^{j_1}
       + \be(\beta) x^{i_2} y^{j_2}
       + \be(\gamma) x^{i_3} y^{j_3}.
$$
Assume further that
all the coefficients of $W$
corresponding to interior lattice points
of the Newton polygon of $W_{e e'}$ vanish.
Then $W_{e e'}$ is the sum of the leading term
and the sub-leading term of $W$
as one puts
$$
 (x, y)
   = (r^{-n(e'')} u, r^{-m(e'')} v),
 \qquad r \in \bR^{>0} \ \text{and} \ 
 u, v \in \bCx,
$$
and take the $r \to \infty$ limit.
Here, $(n(e''), m(e'')) \in M$ is
the primitive outward normal vector
of the edge $e''$ of the Newton polygon of $W_{ee'}$
shown in Figure \ref{fg:asympbone}.
The coamoeba for the sum of three monomials
will be analyzed in Section \ref{sc:coamoeba};
the asymptotic boundaries
coincide with the actual boundaries of the coamoeba,
and the orientations on the asymptotic boundaries
determine which side of the boundary
belongs to the coamoeba.
Hence the orientations of asymptotic boundaries
determine the leading behavior of the coamoeba
near the intersections of asymptotic boundaries
as in Figure \ref{fg:asympbtwo}.

\begin{figure}
\begin{minipage}{.45 \linewidth}
\centering
\input{asymp_b1.pst}
\caption{A pair of adjacent edges of the Newton polygon}
\label{fg:asympbone}
\end{minipage}
\begin{minipage}{.45 \linewidth}
\centering
\input{asymp_b2.pst}
\caption{The leading behavior of the coamoeba
near an intersection of asymptotic boundaries.}
\label{fg:asympbtwo}
\end{minipage}
\end{figure}

\section{Coamoebas for triangles}
 \label{sc:coamoeba}

We discuss the coamoeba of $W^{-1}(0) \subset \bTv$
in this section,
where $W$ is the sum of three Laurent monimials.
The strategy is to reduce to the simplest case
$W(x, y) = 1 + x + y$
by a finite cover of the torus.

\begin{theorem} \label{th:coamoeba_triangle}
Let $\Delta$ be a lattice triangle and
$W$ be the Laurent polynomial
obtained as the sum of monomials corresponding to the vertices of $\Delta$.
Then the coamoeba of $W^{-1}(0)$ has the following description:
\begin{itemize}
 \item
The set $\scA$ of asymptotic boundaries,
equipped with the orientation
coming from the outward normal vectors of $\Delta$,
is an admissible arrangement of oriented lines.
 \item
The coamoeba is the union of colored cells and vertices of $\scA$.
 \item
The restriction of the argument map
to the inverse image of a colored cell
is a diffeomorphism.
It is orientation-preserving if the cell is white,
and reversing if the cell is black.
 \item
The inverse image of a vertex of $\scA$
by the argument map
is homeomorphic to an open interval.
\end{itemize}
\end{theorem}

\begin{proof}
We first consider the simplest case
$$
 \Wone(x, y)
  = 1 + x + y
  \in \bC[\None]
$$
where
$
 \None \cong \bZ^2
$
is a free abelian group of rank two.
It gives a regular map
$$
 \bTvone = \Mone \otimes \bCx \to \bC
$$
where $\Mone = \Hom(\None, \bZ)$ is the dual group of $\None$, and
the coamoeba of $\Wone^{-1}(0)$ is a subset of
$\Tvone = \Mone \otimes \bR / \Mone$.
The zero locus of $\Wone$ is obtained by gluing two disks
$$
 D_1 = \{ (x, y) \in (\bCx)^2 \mid \Im(x) > 0, \, y = - 1 - x \}
$$
and
$$
 D_2 = \{ (x, y) \in (\bCx)^2 \mid \Im(x) < 0, \, y = - 1 - x \}
$$
along three intervals
$$
 I_1 = \{ (x, y) \in (\bR)^2 \mid x < -1 , \, y = - 1 - x \},
$$
$$
 I_2 = \{ (x, y) \in (\bR)^2 \mid -1 < x < 0, \, y = - 1 - x \},
$$
and
$$
 I_3 = \{ (x, y) \in (\bR)^2 \mid x > 0, \, y = - 1 - x \}.
$$
Define
$$
 U_1 = \{ (\theta, \phi) \in \Tvone \mid
  0 < \theta < \frac{1}{2}, \,
  - \frac{1}{2} < \phi < - \frac{1}{2} + \theta \}
$$
and
$$
 U_2 = \{ (\theta, \phi) \in \Tvone \mid
  - \frac{1}{2} < \theta < 0, \,
  \frac{1}{2} + \theta < \phi < \frac{1}{2} \}
$$
as in Figure \ref{fg:Z1_alga}.
Then the argument map gives
an orientation-reversing homeomorphism
from $D_1$ to $U_1$ and
an orientation-preserving homeomorphism
from $D_2$ to $U_2$,
and maps
$I_1$, $I_2$, and $I_3$
to
$(\frac{1}{2}, 0)$,
$(\frac{1}{2}, \frac{1}{2})$, and
$(0, \frac{1}{2})$
respectively.

\begin{figure}
\begin{minipage}{.5 \linewidth}
\centering
\input{Z1_alga.pst}
\caption{the coamoeba for $\Deltaone$}
\label{fg:Z1_alga}
\end{minipage}
\begin{minipage}{.5 \linewidth}
\centering
\input{Z1_x-plane.pst}
\caption{the glued surface}
\label{fg:Z1_surface}
\end{minipage}
\end{figure}

Now we discuss the general case.
Assume that $W$ is given by
$$
 W(x, y) = 1 + x^p y^q + x^r y^s
  \in \bC[N]
$$
for some integers $p$, $q$, $r$ and $s$.
Note that one has a commutative diagram
$$
\begin{psmatrix}[colsep=2]
 \bTv & \bTvone \\
  & \bC
 \psset{nodesep=3pt,arrows=->}
 \ncline{1,1}{1,2}^{\psi \otimes \bCx}
 \ncline{1,1}{2,2}<{W}
 \ncline{1,2}{2,2}>{\Wone}
\end{psmatrix}
$$
where
$
 \psi : M \to \Mone
$
is the linear map represented by the matrix
$
  \begin{pmatrix}
   p & q \\
   r & s
  \end{pmatrix}.
$
It follows from the commutative diagram
$$
\begin{CD}
 \bTv @>{\psi \otimes \bCx}>> \bTvone \\
 @V{\Arg}VV @VV{\Arg}V \\
 \Tv @>{\psi \otimes (\bR / \bZ)}>> \Tvone
\end{CD}
$$
that the coamoeba of $W^{-1}(0)$ is
the pull-back by $\psi \otimes (\bR / \bZ)$
of the coamoeba of $\Wone^{-1}(0)$;
$$
 \Arg(W^{-1}(0))
  = (\psi \otimes (\bR / \bZ))^{-1}(\Arg(\Wone^{-1}(0))).
$$
This reduces the general case
to the the case of $\Wone$ discussed above, and
Theorem \ref{th:coamoeba_triangle} is proved.
\end{proof}

{\bf Acknowledgment}:
We thank Kenji Fukaya
for organizing a workshop in Kinosaki
in March 2006
where this joint project has been initiated.
K.~U. thanks
Alastair Craw for introducing him
to the work of Hanany and Vegh,
and providing helpful explanations and valuable comments.
He also thanks Akira Ishii for collaborations
which improved his understanding of dimer models.
M.~Y. thanks Tohru Eguchi and Akishi Kato for helpful comments.
K.~U. is supported by Grant-in-Aid for Young Scientists (No.18840029).

\bibliographystyle{plain}
\bibliography{bibs}

\def\cprime{$'$}
\begin{thebibliography}{10}

\bibitem{Benvenuti-Pando_Zayas-Tachikawa}
Sergio Benvenuti, Leopoldo~A. Pando~Zayas, and Yuji Tachikawa.
\newblock Triangle anomalies from {E}instein manifolds.
\newblock {\em Adv. Theor. Math. Phys.}, 10(3):395--432, 2006.

\bibitem{Broomhead}
Nathan Broomhead.
\newblock Dimer models and {C}alabi-{Y}au algebras.
\newblock arXiv:0901.4662.

\bibitem{Butti-Zaffaroni_RTD}
Agostino Butti and Alberto Zaffaroni.
\newblock {$R$}-charges from toric diagrams and the equivalence of
  {$a$}-maximization and {$Z$}-minimization.
\newblock {\em J. High Energy Phys.}, (11):019, 42 pp. (electronic), 2005.

\bibitem{Butti-Zaffaroni_TGQGT}
Agostino Butti and Alberto Zaffaroni.
\newblock From toric geometry to quiver gauge theory: the equivalence of
  {$a$}-maximization and {$Z$}-minimization.
\newblock {\em Fortschr. Phys.}, 54(5-6):309--316, 2006.

\bibitem{Davison}
Ben Davison.
\newblock Consistency conditions for brane tilings.
\newblock arXiv:0812.4185.

\bibitem{Duffin}
R.~J. Duffin.
\newblock Potential theory on a rhombic lattice.
\newblock {\em J. Combinatorial Theory}, 5:258--272, 1968.

\bibitem{Feng-He-Kennaway-Vafa}
Bo~Feng, Yang-Hui He, Kristian~D. Kennaway, and Cumrun Vafa.
\newblock Dimer models from mirror symmetry and quivering amoebae.
\newblock {\em Adv. Theor. Math. Phys.}, 12(3):489--545, 2008.

\bibitem{Fowler-Rushbrooke}
R.~H. Fowler and G.~S. Rushbrooke.
\newblock An attempt to extend the statistical theory of perfect solutions.
\newblock {\em Trans. Faraday Soc.}, 33:1272 -- 1294, 1937.

\bibitem{Franco-Hanany-Martelli-Sparks-Vegh-Wecht_GTTGBT}
Sebasti{\'a}n Franco, Amihay Hanany, Dario Martelli, James Sparks, David Vegh,
  and Brian Wecht.
\newblock Gauge theories from toric geometry and brane tilings.
\newblock {\em J. High Energy Phys.}, (1):128, 40 pp. (electronic), 2006.

\bibitem{Franco-Hanany-Vegh-Wecht-Kennaway_BDQGT}
Sebasti{\'a}n Franco, Amihay Hanany, David Vegh, Brian Wecht, and Kristian~D.
  Kennaway.
\newblock Brane dimers and quiver gauge theories.
\newblock {\em J. High Energy Phys.}, (1):096, 48 pp. (electronic), 2006.

\bibitem{Franco-Vegh_MSGTDM}
Sebasti{\'a}n Franco and David Vegh.
\newblock Moduli spaces of gauge theories from dimer models: proof of the
  correspondence.
\newblock {\em J. High Energy Phys.}, (11):054, 26 pp. (electronic), 2006.

\bibitem{Ginzburg_CYA}
Victor Ginzburg.
\newblock {C}alabi-{Y}au algebras.
\newblock math.AG/0612139, 2006.

\bibitem{Gulotta}
Daniel~R. Gulotta.
\newblock Properly ordered dimers, {$R$}-charges, and an efficient inverse
  algorithm.
\newblock {\em J. High Energy Phys.}, (10):014, 31, 2008.

\bibitem{Hanany-Kennaway_DMTD}
Amihay Hanany and Kristian~D. Kennaway.
\newblock Dimer models and toric diagrams.
\newblock hep-th/0503149, 2005.

\bibitem{Hanany-Vegh}
Amihay Hanany and David Vegh.
\newblock Quivers, tilings, branes and rhombi.
\newblock {\em J. High Energy Phys.}, (10):029, 35, 2007.

\bibitem{Ishii-Ueda_09}
Akira Ishii and Kazushi Ueda.
\newblock Dimer models and the special {McKay} correspondence.
\newblock arXiv:0905.0059.

\bibitem{Ishii-Ueda_08}
Akira Ishii and Kazushi Ueda.
\newblock On moduli spaces of quiver representations associated with dimer
  models.
\newblock In {\em Higher dimensional algebraic varieties and vector bundles},
  RIMS K\^oky\^uroku Bessatsu, B9, pages 127--141. Res. Inst. Math. Sci.
  (RIMS), Kyoto, 2008.

\bibitem{Kasteleyn}
P.~W. Kasteleyn.
\newblock Dimer statistics and phase transitions.
\newblock {\em J. Mathematical Phys.}, 4:287--293, 1963.

\bibitem{Kato_ZFSFT}
Akishi Kato.
\newblock Zonotopes and four-dimensional superconformal field theories.
\newblock {\em J. High Energy Phys.}, (6):037, 30 pp. (electronic), 2007.

\bibitem{Kenyon_IDM}
Richard Kenyon.
\newblock An introduction to the dimer model.
\newblock In {\em School and {C}onference on {P}robability {T}heory}, ICTP
  Lect. Notes, XVII, pages 267--304 (electronic). Abdus Salam Int. Cent.
  Theoret. Phys., Trieste, 2004.

\bibitem{Kenyon-Okounkov-Sheffield}
Richard Kenyon, Andrei Okounkov, and Scott Sheffield.
\newblock Dimers and amoebae.
\newblock {\em Ann. of Math. (2)}, 163(3):1019--1056, 2006.

\bibitem{Kenyon-Schlenker}
Richard Kenyon and Jean-Marc Schlenker.
\newblock Rhombic embeddings of planar quad-graphs.
\newblock {\em Trans. Amer. Math. Soc.}, 357(9):3443--3458 (electronic), 2005.

\bibitem{Lee-Rey}
Sangmin Lee and Soo-Jong Rey.
\newblock Comments on anomalies and charges of toric-quiver duals.
\newblock {\em J. High Energy Phys.}, (3):068, 21 pp. (electronic), 2006.

\bibitem{Martelli-Sparks-Yau_GDAM}
Dario Martelli, James Sparks, and Shing-Tung Yau.
\newblock The geometric dual of {$a$}-maximisation for toric
  {S}asaki-{E}instein manifolds.
\newblock {\em Comm. Math. Phys.}, 268(1):39--65, 2006.

\bibitem{Martelli-Sparks-Yau_SEMVM}
Dario Martelli, James Sparks, and Shing-Tung Yau.
\newblock Sasaki-{E}instein manifolds and volume minimisation.
\newblock {\em Comm. Math. Phys.}, 280(3):611--673, 2008.

\bibitem{Mercat_DRSIM}
Christian Mercat.
\newblock Discrete {R}iemann surfaces and the {I}sing model.
\newblock {\em Comm. Math. Phys.}, 218(1):177--216, 2001.

\bibitem{Mozgovoy-Reineke}
Sergey Mozgovoy and Markus Reineke.
\newblock On the noncommutative {D}onaldson-{T}homas invariants arising from
  brane tilings.
\newblock arXiv:0809.0117.

\bibitem{Nakamura_HSAGO}
Iku Nakamura.
\newblock Hilbert schemes of abelian group orbits.
\newblock {\em J. Algebraic Geom.}, 10(4):757--779, 2001.

\bibitem{Okounkov-Reshetikhin-Vafa_QCYCC}
Andrei Okounkov, Nikolai Reshetikhin, and Cumrun Vafa.
\newblock Quantum {C}alabi-{Y}au and classical crystals.
\newblock In {\em The unity of mathematics}, volume 244 of {\em Progr. Math.},
  pages 597--618. Birkh\"auser Boston, Boston, MA, 2006.

\bibitem{Ooguri-Yamazaki_ECG}
Hirosi Ooguri and Masahito Yamazaki.
\newblock Emergent {C}alabi-{Y}au geometry.
\newblock {\em Phys. Rev. Lett.}, 102(16):161601, 4, 2009.

\bibitem{Reid_MC}
Miles Reid.
\newblock {M}c{K}ay correspondence.
\newblock alg-geom/9702016.

\bibitem{Stienstra_CPA}
Jan Stienstra.
\newblock Computation of principal {A}-determinants through dimer dynamics.
\newblock arXiv:0901.3681.

\bibitem{Stienstra_dessins}
Jan Stienstra.
\newblock Hypergeometric systems in two variables, quivers, dimers and dessins
  d'enfants.
\newblock In {\em Modular forms and string duality}, volume~54 of {\em Fields
  Inst. Commun.}, pages 125--161. Amer. Math. Soc., Providence, RI, 2008.

\bibitem{Ueda-Yamazaki_toricdP}
Kazushi Ueda and Masahito Yamazaki.
\newblock Homological mirror symmetry for toric orbifolds of toric del {P}ezzo
  surfaces.
\newblock math.AG/0703267.

\end{thebibliography}

\noindent
Kazushi Ueda

Department of Mathematics,
Graduate School of Science,
Osaka University,
Machikaneyama 1-1,
Toyonaka,
Osaka,
560-0043,
Japan.

{\em e-mail address}\ : \  kazushi@math.sci.osaka-u.ac.jp

\ \\

\noindent
Masahito Yamazaki

Department of Physics,
Graduate School of Science,
University of Tokyo,
Hongo 7-3-1,
Bunkyo-ku,
Tokyo,
113-0033,
Japan

{\em e-mail address}\ : \  yamazaki@hep-th.phys.s.u-tokyo.ac.jp

\end{document}